\documentclass[12pt]{article}
\usepackage{amsmath,amssymb}
\def\R{\mathbb{R}}
\def\C{\mathbb{C}}

\def\N{\mathbb{N}}
\def\L{\mathcal{L}}
\def\A{\mathbf{A}}
\def\B{\mathbf{B}}
\def\a{\mathbf{a}}
\def\b{\mathbf{b}}
\def\Be{\mathcal{B}}
\def\Rr{\mathcal{R}}
\newtheorem{theorem}{\hspace*{\parindent}Theorem}
\newtheorem{lemma}{\hspace*{\parindent}Lemma}
\newtheorem{corollary}{\hspace*{\parindent}Corollary}[theorem]

\newcounter{theremark}

\title{Some new facts concerning the delta neutral H function of Fox}
\author{D.B.\:Karp$^{\rm a,b}$\footnote{Corresponding author. E-mail: D. Karp -- \emph{dimkrp@gmail.com}, E.\:Prilepkina --  \emph{pril-elena@yandex.ru}}~~and E.G.\:Prilepkina$^{\rm a,b}$
\\[10pt]\small{\textit{$\phantom{1}^a$Far Eastern Federal University, 8 Sukhanova street, Vladivostok, 690950, Russia}}\\\small{\textit{$\phantom{1}^b$Institute of Applied Mathematics, FEBRAS, 7 Radio Street, Vladivostok,  690041, Russia}}}
\date{}
\begin{document}
\maketitle
\begin{center}
\parbox{12cm}{
\small\textbf{Abstract.}  In this paper we find several new properties of a class of Fox's $H$ functions
which we call delta neutral. In particular, we find an expansion in the neighborhood of finite nonzero singularity and  give
new Mellin transform formulas under a special restriction on parameters. The last result is applied to prove a conjecture
regarding the representing measure for gamma ratio in Bernstein's theorem. Further, we find the weak limit of
measures expressed in terms of the $H$ function which furnishes a regularization method for integrals
containing the delta neutral zero-balanced function of Fox.  We apply this result to extend a recently discovered
integral equation to zero-balanced case. In the last section of the paper we consider a reduced form of this integral equation
for Meijer's $G$ function. This leads to certain expansions believed to be new even in the case
of the Gauss hypergeometric  function.}
\end{center}

\bigskip


\bigskip


\bigskip

\paragraph{1. Introduction and preliminaries.} In his milestone paper \cite{Norlund}  N{\o}rlund constructed the fundamental system of solutions
of the generalized hypergeometric differential equation  ($D=z\frac{d}{dz}$)
\begin{equation}\label{eq:hyper-equation}
\left\{(D-a_1)(D-a_2)\cdots(D-a_p)-z(D+1-b_1)(D+1-b_2)\cdots(D+1-b_p)\right\}y=0
\end{equation}
in the neighborhood of its three regular singular points $0$, $1$, $\infty$ in terms
of the Mellin-Barnes integrals \cite{ParKam}.  Solutions in the neighborhood of $0$, $\infty$ were previously built by Thomae in 1870
in terms of the generalized hypergeometric functions. Rather straightforward examination shows that  N{\o}rlund's solution in the neighborhood
of $1$ that corresponds to the local exponent $\psi-1=\sum_{k=1}^{p}(b_k-a_k)-1$ is equal to Meijer's $G$ function $G^{p,0}_{p,p}$ defined
by the contour integral (\ref{eq:Fox}) below on setting $\A=\B=(1,\ldots,1)$.  See also \cite[section~8.2]{PBM3} for its definition and
properties.  In his paper  N{\o}rlund discovered a large number of expansions, representations
and connection formulas for his solutions. Among many identities for the function $G^{p,0}_{p,p}$ found in \cite{Norlund}
a prominent role is played by the following expansion \cite[(1.33)]{Norlund}:
\begin{equation}\label{eq:Norlund1}
G^{p,0}_{p,p}\!\left(\!z~\vline\begin{array}{l}\b\\\a\end{array}\!\!\right)=\frac{z^{a_k}(1-z)^{\psi-1}}{\Gamma(\psi)}
\sum\limits_{n=0}^{\infty}\frac{g_n(\a_{[k]};\b)}{(\psi)_n}(1-z)^n,~~~k=1,2,\ldots,p.
\end{equation}
The coefficients $g_n(\a_{[k]};\b)$ satisfy two different recurrence relations (in $p$ and $n$) and the following explicit
formula \cite[(1.28), (2.7), (2.11)]{Norlund}:
\begin{equation}\label{eq:Norlund-explicit}
g_n(\a_{[p]};\b)=\sum\limits_{0\leq{j_{1}}\leq{j_{2}}\leq\cdots\leq{j_{p-2}}\leq{n}}
\prod\limits_{m=1}^{p-1}\frac{(\psi_m+j_{m-1})_{j_{m}-j_{m-1}}}{(j_{m}-j_{m-1})!}(b_{m+1}-a_{m})_{j_{m}-j_{m-1}},
\end{equation}
where $\psi_m=\sum_{i=1}^{m}(b_i-a_i)$, $j_0=0$, $j_{p-1}=n$, and
$$
\a=(a_1,\ldots,a_p),~~\b=(b_1,\ldots,b_p),~~\a_{[k]}=(a_1,\ldots,a_{k-1},a_{k+1},\ldots,a_p).
$$
The coefficient $g_n(\a_{[k]};\b)$ is obtained from $g_n(\a_{[p]};\b)$ by exchanging the roles of $a_p$ and $a_k$.
N{\o}rlund  also mentioned without reference or proof that ``as $z=1$ is a regular singularity, (\ref{eq:Norlund1})
is convergent in the circle $|z-1|<1$''.  Indeed, this fact follows from the general theory of Fuchsian differential
equations, see, for instance, \cite[Theorem~11.3]{Hartman}. Convergence of (\ref{eq:Norlund1}) in the disk $|z-1|<1$ allows using this series
for analytic continuation of $G^{p,0}_{p,p}$ from inside of the unit disk to the disk $|z-1|<1$.  One of the consequences of (\ref{eq:Norlund1})
found in \cite{Norlund} is the following little-known formula for the Mellin transform of the $G$ function valid
for $\psi=\sum_{j=1}^{p}(b_j-a_j)=-m$, $m\in\N_0=\N\cup\{0\}$:
\begin{equation}\label{eq:GMellinNorlund}
\int\limits_{0}^{1}x^{s-1}G^{p,0}_{p,p}\!\left(\!x~\vline\begin{array}{l}\b\\\a\end{array}\!\!\right)dx
=\prod\limits_{j=1}^{p}\frac{\Gamma(a_j+s)}{\Gamma(b_j+s)}-q(s),~~~\Re(s)>-\Re(a_j),~j=1,\ldots,p.
\end{equation}
Here $q(s)$ is the polynomial of degree $m$ given by
$$
q(s)=\sum\limits_{j=0}^{m}g_{m-j}(\a_{[k]};\b)(s+a_k-j)_j,~~~k=1,2,\ldots,p,
$$
where the standard notation $(\alpha)_j=\Gamma(\alpha+j)/\Gamma(\alpha)$ for the rising factorial has been used.
Note that for $\psi>0$ the term $q(s)$ is absent from (\ref{eq:GMellinNorlund}) while for $\psi\le0$,
$-\psi\notin\N_0$ the integral on the left diverges.  An account of  N{\o}rlund's work along with
extensions and applications can be found in our recent paper \cite{KPNorlund}, where an interested reader is referred to
for details.

The first motivation for this paper is to find generalizations of formulas (\ref{eq:Norlund1}) and (\ref{eq:GMellinNorlund})
to a particular case of Fox's $H$ function defined by
\begin{equation}\label{eq:Fox}
H_{q,p}^{p,0}\left(z\left|\begin{array}{l} (\B,\b)\\(\A,\a)\end{array}\right.\right)=
\frac{1}{2\pi{i}}\int\limits_{\L}\frac{\prod\nolimits_{k=1}^{p}\Gamma(A_ks+a_k)}
{\prod\nolimits_{j=1}^{q}\Gamma(B_j s+b_j)} z^{-s}ds.
\end{equation}
Here $\A=(A_1,\ldots,A_p),$ $\B=(B_1,\ldots,B_q)$ are positive vectors while
$\a=(a_1,\ldots,a_p),$ $\b=(b_1,\ldots,b_q)$ are arbitrary real vectors
(most results of this paper remain valid for complex $\a$, $\b$ by analytic continuation, but we restrict
our attention to real case for simplicity).  The contour $\L$  starts and ends at infinity and leaves all
the poles of the integrand on the left.  It can be one of the following:
\begin{itemize}

\item $\L=\L_{-\infty}$ is a left loop situated in a horizontal strip starting at the point
$-\infty+i\varphi_1$ and terminating at the point $-\infty+i\varphi_2$, where  $-\infty<\varphi_1<\varphi_2<\infty$;

\item $\L=\L_{\infty}$ is a right loop situated in a horizontal strip starting at the point
 $\infty+i\varphi_1$  and terminating at the point $\infty+i\varphi_2$, where $-\infty<\varphi_1<\varphi_2<\infty$;

\item $\L=\L_{ic}$ is the vertical line $\Re{z}=c$, where $c>\gamma=-\min\limits_{1\leq k\leq p}({a_k}/{A_k})$.

\end{itemize}
Details regarding the contour and conditions for convergence of the integral in (\ref{eq:Fox})
can be found in \cite[sections 1.1,1.2]{KilSaig}. In this paper we deal with the case of the $H$ function
such that all the poles of the integrand lie in some left half-plane which permits considering
$\L_{\pm\infty}$ as the boundary of a convex simply connected domain.

One of the key parameters characterizing the $H$ function is the difference
$\Delta=\sum_{j=1}^{q}B_j-\sum_{i=1}^{p}A_i$.  Obviously, $\Delta=0$ for $G_{p,p}^{p,0}$.  Moreover, only
for $\Delta=0$  has the function  $H_{q,p}^{p,0}$ a compact support as does $G_{p,p}^{p,0}$, see \cite{KPCMFT}.
The same paper contains further arguments why $\Delta=0$ is necessary for $H_{q,p}^{p,0}$ to be analogous to $G_{p,p}^{p,0}$.
Hence, in what follows  we will assume that
\begin{equation}\label{eq:balance}
\Delta=\sum\nolimits_{j=1}^{q}B_j-\sum\nolimits_{i=1}^{p}A_i=0.
\end{equation}
Under this condition the function $H_{q,p}^{p,0}$ will be called delta neutral (we avoid here the term ``balanced'' since it is
reserved for the case $\mu=0$, where $\mu$ is defined in (\ref{eq:mu}) below).  It follows from  \cite[Theorems~3.3]{KilSaig}
and \cite[Theorem~6]{KPCMFT} that the delta neutral $H_{q,p}^{p,0}(z)$ is supported on the disk $|z|\leq\rho$, where the number
\begin{equation}\label{eq:rho}
\rho=\prod\limits_{k=1}^{p}A_k^{A_k}\prod\limits_{j=1}^{q}B_j^{-B_j}
\end{equation}
plays the same role for the delta neutral $H_{q,p}^{p,0}$ unity plays for $G_{p,p}^{p,0}$. Further, define
\begin{equation}\label{eq:mu}
\mu=\sum\nolimits_{j=1}^{q}b_j-\sum\nolimits_{k=1}^{p}a_k+\frac{p-q}{2}.
\end{equation}
This parameter generalizes the parameter $\psi$ of the $G$ function.  Suppose $0<x<\rho$. Then choosing
the contour $\L=\L_{-\infty}$ in (\ref{eq:Fox}) guarantees convergence of the integral for all complex $\mu$; if $\mu>0$
the integral over $\L=\L_{ic}$ also converges  and has the same value.  If $x>\rho$ then choosing $\L=\L_{+\infty}$
in (\ref{eq:Fox}) guarantees convergence of the integral for all complex $\mu$; if $\mu>0$ the integral over $\L=\L_{ic}$
also converges  and has the same value.  See \cite[Theorems~1.1~and~3.3]{KilSaig} and \cite[Theorem~6]{KPCMFT} for the
proof of these claims.

Next motivation for this paper comes from the fact that the delta neutral $H$ function and its particular case $G^{p,0}_{p,p}$
play a noticeable role in statistics and have been used there since 1932 paper of Wilks \cite{Wilks} which appeared
before Meijer and Fox introduced their functions.  In that paper Wilks introduced an integral equation which he
called type B integral equation. Many likelihood ratio criteria in multivariate hypothesis testing have probability density
functions satisfying this integral equation including general linear hypothesis, sphericity hypothesis,
independence of groups of variables, equality of several covariance matrices and other.
See details and numerous reference in \cite{And,GuptTang}.  A quick examination reveals that the solution of Wilks'
type B integral equation is precisely $G^{p,0}_{p,p}$.  Seventeen years after Wilks paper, Box \cite{Box} considered
even more general integral equation whose solution is easily seen to be the delta neutral $H$ function
(introduced by Fox in \cite{Fox} 12 years after Box's paper). The role this equation plays in multivariate statistics is explained in
section~8.5 of Anderson's authoritative volume \cite{And}. In a series of papers Gupta and Tang
tried to solve both Wilks' and Box's equations in terms of series expansions and presented further applications of their solutions
in statistics. See \cite{GuptTang,TangGupta1,TangGupta2} and numerous references therein. In particular, in \cite{GuptTang} they
found two series  expansions for ``a statistic with general moment function'' ($=$ the delta neutral $H$ function).
Their second expansion corresponds to our Theorem~\ref{th:NorforH} below, but without convergence proof and
with more complicated formulas for the coefficients.  They also base their argument on an incorrect theorem
of Nair \cite[Lemma~4.1]{GuptTang}.  Furthermore, the delta neutral $H$ function represents the probability density function of a
product of positive powers of beta distributed random variables, see \cite[section~4.2]{MSH}. Chapter~4 of this book also
contains further applications of the $H$ function in statistics.  This line of research continues today -
see, for instance, \cite{CoelhoArnold} and references there.  Independently, probability distribution with  the $G$ function density
has been found to be the stationary distribution of certain Markov chains and is known as Dufresne law in this context. See
\cite{ChamLetac,Dufresne1998,Dufresne2010} and references therein.

The first purpose of this paper is to give a complete description of the singularity of $H_{q,p}^{p,0}(z)$ at the point $z=\rho$.
We establish an expansion in the neighborhood of this point which generalizes (\ref{eq:Norlund1}).
The main difference from the $G$ function case lies in the fact that the function $H_{q,p}^{p,0}(z)$ does not satisfy
any reasonable differential equation, so that no Frobenius method is available for computation of the coefficients
and no general theory can be resorted to for convergence proof.  We overcome both of these difficulties.  We also give a direct independent proof of
the Mellin transform formula generalizing (\ref{eq:GMellinNorlund}). The role of $\psi$ is thereby played by the parameter $\mu$
from (\ref{eq:mu}).

Recall that a function $f:(0,\infty)\to(0,\infty)$ is called completely
monotonic if $(-1)^nf^{(n)}(x)\ge0$ for $x>0$ and $n\in\N_0$. The celebrated Bernstein theorem asserts that completely monotonic
functions are precisely those that can be expressed by the Laplace transform of a nonnegative measure.
In our recent paper  \cite{KPCMFT} we found conditions for the ratio
$$
W(x)=\frac{\prod\nolimits_{k=1}^{p}\Gamma(A_kx+a_k)}{\prod\nolimits_{j=1}^{q}\Gamma(B_jx+b_j)}
$$
to be completely monotonic on $(0,\infty)$. The second motivation for this paper comes from the fact that the delta neutral $H_{q,p}^{p,0}(z)$ constitutes the representing
measure for $W(x)$ in Bernstein theorem if $\mu>0$, i.e. \cite[Theorem~7]{KPCMFT}
\begin{equation}\label{eq:measureH}
W(x)=\int_{\log(1/\rho)}^{\infty}e^{-tx}H_{q,p}^{p,0}\left(e^{-t}\left|\begin{array}{l}(\B,\b)\\(\A,\a)\end{array}\right.\right)\!dt.
\end{equation}
We conjectured in \cite[Conjecture~2]{KPCMFT} that for $\mu=0$ the representing measure takes the form
\begin{equation}\label{eq:limitmeasure}
A\delta_{\log(1/\rho)}+H_{q,p}^{p,0}\left(e^{-t}\left|\begin{array}{l}(\B,\b)\\(\A,\a)\end{array}\right.\right)dt,
\end{equation}
where $\delta_{x}$ denotes the unit mass concentrated at $x$. In this paper we prove this conjecture and find an explicit
expression for the constant $A$.  Moreover, we go one step further and show that the measure in (\ref{eq:measureH}) converges
weakly to the measure in (\ref{eq:limitmeasure}) as $\mu\to0$ even when $W(x)$ is not completely monotonic and both measures are signed.
This provides a regularization method  for the integrals containing the delta neutral function $H_{q,p}^{p,0}$ which is also zero-balanced, i.e. for $\mu=0$.
In Corollary~\ref{cr:intequation} below we apply this regularization to extend the integral equation
\begin{equation}\label{eq:ident1}
H_{q,p}^{p,0}\left(x\left|\begin{array}{l}(\B,\b)\\(\A,\a)\end{array}\right.\right)=\frac{1}{\log(\rho/x)}\int_{x/\rho}^{1}
H_{q,p}^{p,0}\left(\frac{x}{u}\left|\begin{array}{l}(\B,\b)\\(\A,\a)\end{array}\right.\right)\frac{Q(u)}{u}du,
\end{equation}
to the zero-balanced case.  This integral equation was established in \cite[Theorem~8]{KPCMFT} for $x\in(0,\rho)$ under the
restrictions $\A,\B>0$, $\a,\b\ge0$, $\rho\le1$ and $\mu>0$.
Here
\begin{equation}\label{eq:Qu-def}
Q(u)=\sum\limits_{i=1}^{p}\frac{u^{a_i/A_i}}{1-u^{1/A_i}}-\sum\limits_{j=1}^{q}\frac{u^{b_j/B_j}}{1-u^{1/B_j}}~~
\text{for}~~u\in(0,1).
\end{equation}
When $H_{q,p}^{p,0}$ reduces to Meijer's $G$ function, the integral on the right hand side of (\ref{eq:ident1}) can be evaluated.
This leads to a new functional-differential relation and a new expansion for Meijer's $G$ function $G_{p,p}^{p,0}$
which we present in the ultimate section of this paper.

\paragraph{2. Main results for the $H$ function.}
Let us start by recalling that $n$-th Bernoulli polynomial $\Be_{n}(a)$ is defined via the generating function
\cite[p.588]{NIST}
$$
\frac{te^{at}}{e^t-1}=\sum\limits_{n=0}^{\infty}\Be_n(a)\frac{t^n}{n!}, ~~~~|t|<2\pi.
$$

The next lemma contains a corrected form of the asymptotic formula found in \cite{GuptTang}
which was originally derived in \cite{Box} in a slightly different form (see also \cite[8.5.1]{And}).
This formula follows immediately from Hermite-Barnes asymptotic expansion for the gamma function, see, for instance,
\cite[(1.6)]{Nemes}.  Details can be found in our forthcoming paper \cite{KNPGammaPaper}.
\begin{lemma}\label{lm:GammaPoincareExp}
Suppose $\Delta=0$, $M\in\N_0$ and $a_k,b_k$ are arbitrary real numbers. Then
\begin{equation}\label{eq:asimptotic}
\frac{\prod\nolimits_{k=1}^{p}\Gamma(A_kz+a_k)}{\prod\nolimits_{j=1}^{q}\Gamma(B_jz+b_j)}
=\nu\rho^{z}z^{-\mu}\left(\sum\limits_{r=0}^M\frac{l_r}{z^r}+O\left(\frac{1}{z^{M+1}}\right)\right)~\text{as}~z\to\infty
\end{equation}
in the sector $|\arg{z}|<\pi-\delta$, $0<\delta<\pi$.  Here
\begin{equation}\label{eq:nu}
\nu=(2\pi)^{(p-q)/2}\prod\nolimits_{k=1}^{p}A_k^{a_k-1/2}\prod\nolimits_{j=1}^{q}B_k^{1/2-b_j},
\end{equation}
$\rho$ and  $\mu$ are defined in \emph{(\ref{eq:rho})} and  \emph{(\ref{eq:mu})}, respectively.
The coefficients $l_r$ satisfy the recurrence relation
\begin{equation}\label{eq:lr}
l_r=\frac{1}{r}\sum\limits_{m=1}^r q_m l_{r-m},~~\text{with}~~l_0=1
\end{equation}
and
\begin{equation}\label{eq:qt}
q_m=\frac{(-1)^{m+1}}{m+1}\left[\sum\limits_{k=1}^p\frac{\Be_{m+1}(a_k)}{A_k^m}-\sum\limits_{j=1}^q\frac{\Be_{m+1}(b_j)}{B_j^m}\right].
\end{equation}
\end{lemma}
\textbf{Remark.} It is known \cite[Lemma~1]{Kalinin} that the recurrence (\ref{eq:lr}) can be solved to give the following explicit expressions for $l_r$:
$$
l_r=\sum\limits_{k_1+2k_2+\cdots+rk_r=r}\frac{q_1^{k_1}(q_2/2)^{k_2}\cdots (q_r/r)^{k_r}}{k_1!k_2!\cdots k_r!}
=\sum\limits_{n=1}^{r}\frac{1}{n!}\sum\limits_{k_1+k_2+\cdots+k_n=r}\prod\limits_{i=1}^{n}\frac{q_{k_i}}{k_i}.
$$
Moreover, Nair \cite[section~8]{Nair} found a determinantal expression
for such solution which in our notation takes the form
$$
l_r=\frac{\det(\Omega_r)}{r!},~~~\Omega_r=[\omega_{i,j}]_{i,j=1}^{r},~~\omega_{i,j}\!=\!\left\{\!\!\begin{array}{ll}q_{i-j+1}(i-1)!/(j-1)!, & i\ge{j},\\
-1, &i=j-1,\\0, &i<j-1.
\end{array}\right.
$$

Will need the so called signless non-central Stirling numbers of the first kind $s_{\sigma}(n,l)$ \cite[8.5]{Charalambides} defined by
their their ''horizontal'' generating function \cite[(A.2)]{Weniger2010} as follows:
\begin{equation*}
(\sigma+x)_{n}=\sum\limits_{l=0}^{n}x^{l}s_{\sigma}(n,l),~~\text{where}~(a)_n=\Gamma(a+n)/\Gamma(a).
\end{equation*}
These numbers have been studied by Carlitz \cite{Carlitz1,Carlitz2}.
In particular, Carlitz found another nice ''double'' generating function \cite[(5.4)]{Carlitz1}
$$
\sum\limits_{l,n=0}^{\infty}s_{\sigma}(n,l)y^l\frac{x^n}{n!}=(1-x)^{-\sigma-y}
$$
and used it to relate $s_{\sigma}(n,l)$ to Bernoulli-N{\o}rlund (or generalized Bernoulli) polynomials $\Be^{(\sigma)}_{k}(x)$
 defined  by the generating function \cite[(1)]{Norlund61}
$$
\frac{t^{\sigma}e^{xt}}{(e^t-1)^{\sigma}}=\sum\limits_{k=0}^{\infty}\Be^{(\sigma)}_{k}(x)\frac{t^k}{k!}.
$$
His expression as given in \cite[(7.6)]{Carlitz2} is:
$$
s_{\sigma}(n,l)=\binom{-l-1}{n-l}\Be^{(n+1)}_{n-l}(1-\sigma)=\frac{(-1)^{n-l}(l+1)_{n-l}}{(n-l)!}\Be^{(n+1)}_{n-l}(1-\sigma).
$$
Numerous formulas  for these numbers are also collected in \cite{Charalambides}.  The key role in the proof of our first
theorem will be played by the next statement which is a combination of \cite[Theorems~III and IV]{Norlund14} (see also a related result in
\cite[\S94~I,II]{Nielsen}).

\renewcommand{\thetheorem}{\Alph{theorem}}
\setcounter{theorem}{0}

\begin{theorem}\label{th:Norlund-Nielsen}
Let $b_k,\beta\in\C$. Suppose the series
$$
\Omega(z)=\sum\limits_{k=0}^{\infty}\frac{b_{k}k!}{(z+\beta)_{k+1}}
$$
converges for some finite $z_0\in\C$.  Then it converges uniformly  in the half-plane $\Re(z-z_0)\ge\varepsilon$
for any $\varepsilon>0$ and absolutely in the half-plane $\Re(z-z_0)>1$.
 The function $\Omega(z)$ admits the representation
$$
\Omega(z)=\int_0^1 t^{z+\beta-1}\phi(t)dt,~~\text{where}~~\phi(t)=\sum_{k=0}^{\infty}b_k(1-t)^k
$$
is analytic in the disk $|1-t|<1$.
\end{theorem}

\renewcommand{\thetheorem}{\arabic{theorem}}
\setcounter{theorem}{0}

The following theorem extends N{\o}rlund's expansion (\ref{eq:Norlund1}) to the  delta neutral $H$ function.
\begin{theorem}\label{th:NorforH}
Suppose $\Delta=0$.  For arbitrary real $\theta$ and complex $t$ lying in the intersection  of the disks
$|t|<1$ and $|t-1|<1$ the following representations hold true\emph{:}
\begin{align}\label{eq:Hexpansion1}
H_{q,p}^{p,0}\left(\rho{t}\left|\begin{array}{l}(\B,\b)\\(\A,\a)\end{array}\right.\right)
=&{\nu}t^{\theta+1}(1-t)^{\mu-1}\sum\limits_{n=0}^{\infty}\frac{(1-t)^{n}}{\Gamma(\mu+n)}\sum\limits_{r=0}^{n}c_{r}s_{\theta+\mu}(n,r)
\\\label{eq:Hexpansion2}
={\nu}t^{\theta+1}(1-t)^{\mu-1}\sum\limits_{n=0}^{\infty}(1-t)^{n}
&\sum\limits_{r+k=n}\frac{(-1)^kl_r}{k!\Gamma(r+\mu)}\Be^{(n+\mu)}_{k}(-\theta),
\end{align}
where $\nu$ is given in \emph{(\ref{eq:nu})}, $l_r$ is defined in \emph{(\ref{eq:lr})}, $c_0=1$ and $c_r$ for $r\ge1$
are found from the recurrence
\begin{equation}\label{eq:ñr}
c_r=\frac{1}{r}\sum\limits_{m=1}^r \widetilde{q}_m c_{r-m},
\end{equation}
\begin{equation}
\widetilde{q}_m= \frac{(-1)^{m+1}}{m+1}
\left[\Be_{m+1}(\theta+\mu)-{\Be}_{m+1}(\theta+1)+\sum\limits_{k=1}^p\frac{\Be_{m+1}(a_k)}{A_k^m}
-\sum\limits_{j=1}^q\frac{\Be_{m+1}(b_j)}{B_j^m}\right].
\end{equation}
The series in \emph{(\ref{eq:Hexpansion1}), (\ref{eq:Hexpansion2})} converge in the disk $|t-1|<1$ and their sum
represents the analytic continuation of the left hand side of \emph{(\ref{eq:Hexpansion1})} to the disk $|t-1|<1$.
\end{theorem}

\textbf{Remark.}  Expansion (\ref{eq:Hexpansion1}) shows that the delta neutral $H$ function on
the left of (\ref{eq:Hexpansion1}) has a branch point at $t=1$ for all non-integer values of $\mu$ (algebraic or
logarithmic depending on rationality of $\mu$).  For all integer $\mu$, including negative values, this function is analytic
around $t=1$.

\textbf{Proof.}  First, we prove expansion (\ref{eq:Hexpansion1}).
In our forthcoming paper \cite{KNPGammaPaper} we establish the formula
\begin{equation}\label{eq:GammaRatioExt}
\rho^{-z}\frac{\prod\nolimits_{k=1}^{p}\Gamma(A_kz+a_k)}{\prod\nolimits_{j=1}^{q}\Gamma(B_jz+b_j)}
\frac{\Gamma(z+\theta+\mu)}{\Gamma(z+\theta+1)}
=\sum\limits_{n=0}^{\infty}\frac{n!a_n}{(z+\theta+\mu)_{n+1}},
\end{equation}
where
\begin{equation}\label{eq:GammaRatioExt1}
a_n =\frac{\nu}{n!}\sum\limits_{r=0}^{n}c_{r}s_{\theta+\mu}(n,r)
\end{equation}
with $\nu$ defined in (\ref{eq:nu}), and the series in (\ref{eq:GammaRatioExt}) converges in some right half-plane $\Re{z}>\lambda$.
Therefore, according to Theorem~\ref{th:Norlund-Nielsen} the series $\sum_{n=0}^{\infty}a_n(1-t)^n$
converges in the disk $|1-t|<1$. This, in turn, implies that the series
$$
\sum_{n=0}^{\infty}\frac{a_nn!}{\Gamma(\mu+n)}(1-t)^n
$$
converges in the same disk (because $n!/\Gamma(\mu+n)\sim{n^{1-\mu}}$ as $n\to\infty$).
Next, we assume temporarily $\mu>0$ and $\Re{z}>\lambda+1$, multiply the above series by $t^{z+\theta}(1-t)^{\mu-1}$ and integrate
term by term using (\ref{eq:GammaRatioExt}):
\begin{multline*}
\int_0^1\left\{\sum_{n=0}^{\infty}\frac{a_nn!}{\Gamma(\mu+n)}t^{z+\theta}(1-t)^{n+\mu-1}\right\}dt
=\sum_{n=0}^{\infty}\frac{a_nn!}{\Gamma(\mu+n)}\int_0^1t^{z+\theta}(1-t)^{n+\mu-1}dt
\\
=\sum_{n=0}^{\infty}\frac{a_nn!\Gamma(z+\theta+1)}{\Gamma(z+\theta+\mu+n+1)}
=\rho^{-z}\frac{\prod\nolimits_{k=1}^{p}\Gamma(A_kz+a_k)}{\prod\nolimits_{j=1}^{q}\Gamma(B_jz+b_j)}
=\int_{0}^{1}t^{z-1}H_{q,p}^{p,0}\left(\rho{t}\left|\begin{array}{l}(\B,\b)\\(\A,\a)\end{array}\right.\right)dt,
\end{multline*}
where the last equality is the standard formula for the Mellin transform of the $H$ function \cite[Theorem~2.2]{KilSaig} combined with the fact that the delta neutral $H$
function has compact support, see \cite[Theorem~6]{KPCMFT}. In view of uniqueness of the inverse Mellin transform comparison of the first and the last term and
formula (\ref{eq:GammaRatioExt1}) yield expansion (\ref{eq:Hexpansion1}).  It remains to justify the termwise integration above.   According to a version of Lebesgue
dominated convergence theorem  \cite[Theorem~3.6.2]{Bendetto} it suffices to demonstrate that
$$
\sum_{n=0}^{\infty}\frac{|a_n|n!}{\Gamma(\mu+n)}\int_0^1|t^{z+\theta}(1-t)^{n+\mu-1}|dt<\infty.
$$
But this series is clearly equal to
$$
\sum_{n=0}^{\infty}\frac{|a_n|n!\Gamma(x+\theta+1)}{\Gamma(x+\theta+\mu+n+1)}=\frac{\Gamma(x+\theta+1)}{\Gamma(x+\theta+\mu)}\sum_{n=0}^{\infty}\frac{|a_n|n!}{(x+\theta+\mu)_{n+1}},
$$
where $x=\Re{z}$. According to Theorem~\ref{th:Norlund-Nielsen} it converges for $\Re{z}>\lambda+1$.
This proves the validity of expansion (\ref{eq:Hexpansion1}) for $\mu>0$. Recalling the definition of $\mu$ in (\ref{eq:mu})
we see that both sides of (\ref{eq:Hexpansion1}) are well defined analytic functions of $a_i$, $b_j$.
Hence, the restriction $\mu>0$ can be removed by analytic continuation.

To prove (\ref{eq:Hexpansion2}) substitute the expansion \cite[(43)]{Norlund61}
$$
\frac{1}{z^{\beta}}=\sum\limits_{k=0}^{\infty}\frac{(-1)^k\Be^{(\beta+k)}_{k}(-\theta)(\beta)_{k}\Gamma(z+\theta+1)}{k!\Gamma(z+\theta+\beta+k+1)},
$$
convergent in the half plane $\Re{z}>0$, for the powers $z^{-r-\mu}$ in asymptotic formula (\ref{eq:asimptotic}) and rearrange.
These manipulations yield
\begin{multline*}
\rho^{-z}\frac{\prod\nolimits_{k=1}^{p}\Gamma(A_kz+a_k)}{\prod\nolimits_{j=1}^{q}\Gamma(B_jz+b_j)}
\sim\nu\sum\limits_{r=0}^\infty\frac{l_r}{z^{r+\mu}}=\nu\sum\limits_{r=0}^{\infty}l_r
\sum\limits_{k=0}^{\infty}\frac{(-1)^k\Be^{(r+\mu+k)}_{k}(-\theta)(r+\mu)_{k}\Gamma(z+\theta+1)}{k!\Gamma(z+\theta+r+\mu+k+1)}
\\
=\nu\sum\limits_{n=0}^{\infty}\frac{\Gamma(z+\theta+1)\Gamma(n+\mu)}{\Gamma(z+\theta+\mu+n+1)}\sum\limits_{r+k=n}
\frac{(-1)^kl_r}{k!\Gamma(r+\mu)}\Be^{(n+\mu)}_{k}(-\theta).
\end{multline*}
It is known \cite[III, p.454]{Nielsen1904} and is easy to see that the resulting series is also asymptotic for the left hand side.
On the other hand, it clearly coincides with the series (\ref{eq:GammaRatioExt}). Uniqueness of inverse factorial series
(convergent or asymptotic) shows that we have thus merely obtained an alternative expression for the coefficients $a_n$:
$$
a_n=\frac{\nu}{n!}\Gamma(n+\mu)\sum\limits_{r+k=n}\frac{(-1)^kl_r}{k!\Gamma(r+\mu)}\Be^{(n+\mu)}_{k}(-\theta).
$$
Substituting this expression into (\ref{eq:Hexpansion1}) in place of (\ref{eq:GammaRatioExt1}) we
immediately arrive at (\ref{eq:Hexpansion2}). $\hfill\square$

The next theorem contains a generalization of N{\o}rlund's formula (\ref{eq:GMellinNorlund}).
\begin{theorem}\label{th:Fox}
Suppose $\Delta=0$ and $\mu=-m$, $m\in\N_0$. Let $\gamma=-\min\limits_{1\leq k\leq p}({a_k}/{A_k})$ denote the rightmost pole of the integrand in \emph{(\ref{eq:Fox})}.
Then the Mellin transform of the delta neutral $H$ function exists in the half-plane $\Re{s}>\gamma$ and is given by\emph{:}
$$
\int_{0}^{\rho}H_{q,p}^{p,0}\left(x\left|\begin{array}{l}(\B,\b)\\(\A,\a)\end{array}\right.\right)x^{s-1}dx
=\frac{\prod\nolimits_{k=1}^{p}\Gamma(A_ks+a_k)}{\prod\nolimits_{j=1}^{q}\Gamma(B_js+b_j)}
-\nu\rho^s\sum\limits_{k=0}^m l_{m-k}s^k,
$$
where the coefficients $l_{r}$ are computed by \emph{(\ref{eq:lr})}.
\end{theorem}

\textbf{Remark.} Theorem~\ref{th:NorforH} implies that the Mellin transform of $H$ in
Theorem~\ref{th:Fox} does not exist for any value of $s$  if $\mu$ is negative and non-integer.

\textbf{Proof.} Note first that the contour $\L_{-\infty}$ in the definition of the $H$ function can be placed on the left of the
vertical line $\Re{z}=c$, $c>\gamma$.  Further, the part of the contour outside of some disk $|s|<r_0$ can be deformed into the rays
$\Rr_{\pm\varepsilon}=\{re^{i(\pi\pm\varepsilon)}:~r\ge{r_0}>0\}$ for some $\varepsilon>0$.  To justify this claim we need to demonstrate that
the slanted ray $\Rr_{-\varepsilon}$ can be deformed into the horizontal ray
$\Rr=\{t+ir_0\sin\varepsilon:~-\infty<t\leq-r_0\cos\varepsilon\}$ without altering
the value of the integral in the definition of the  $H$ function. Indeed, using the reflection formula for the gamma function we get
\begin{equation}\label{eq:Gammaratio}
\frac{\prod_{i=1}^{p}\Gamma(A_is+a_i)}{\prod_{j=1}^{q}\Gamma(B_js+b_j)}
=\pi^{p-q}\frac{\prod_{j=1}^{q}\Gamma(1-B_js-b_j)}{\prod_{i=1}^{p}\Gamma(1-A_is-a_i)}
\frac{\prod_{j=1}^{q}\sin(\pi(B_js+b_j))}{\prod_{i=1}^{p}\sin(\pi(A_is+a_i))}.
\end{equation}
In view of  condition (\ref{eq:balance}) it is straightforward to see that the  quotient of the sine functions
above is bounded in the domain $|\Im{s}|\ge{r_0}\sin\varepsilon$ (see, for instance, \cite[proof of theorem~3.3]{KilSaig}).
Asymptotic expansion (\ref{eq:asimptotic}) applied to the function
$$
t\to x^{-s}\frac{\prod_{j=1}^{q}\Gamma(1-B_j(-u+it)-b_j)}{\prod_{i=1}^{p}\Gamma(1-A_i(-u+it)-a_i)}
$$
on the vertical segment $\Rr_u=[-u+ir_0\sin\varepsilon,\, -u+iu\tan\varepsilon]$ connecting $\Rr$ and  $\Rr_{-\varepsilon}$ implies that
$$
\lim\limits_{u\to+\infty} \int\limits_{\Rr_u}\frac{\prod\nolimits_{k=1}^{p}\Gamma(A_ks+a_k)}{\prod\nolimits_{j=1}^{q}\Gamma(B_j s+b_j)} x^{-s}ds=0,
~s=-u+it.
$$
Then, according to the residue theorem
$$
\int\limits_{\Rr_{-\varepsilon}}\frac{\prod\nolimits_{k=1}^{p}\Gamma(A_ks+a_k)}{\prod\nolimits_{j=1}^{q}\Gamma(B_j s+b_j)}x^{-s}ds
=\int\limits_{\Rr}\frac{\prod\nolimits_{k=1}^{p}\Gamma(A_ks+a_k)}{\prod\nolimits_{j=1}^{q}\Gamma(B_j s+b_j)} x^{-s}ds.
$$
A similar argument shows that $\Rr_{+\varepsilon}$ can be deformed into a horizontal line.

Next, set
$$
P_m(s)=\nu  \sum\limits_{k=0}^m l_{m-k}s^k~~\text{and}~~\mathcal{H}(s)=\frac{\prod\nolimits_{k=1}^{p}\Gamma(A_ks+a_k)}
{\prod\nolimits_{j=1}^{q}\Gamma(B_j s+b_j)}-\rho^s P_m(s).
$$
The function $\mathcal{H}(s)$ has no poles in the region bounded by  $\L_{-\infty}$ and  $\L_{ic}$.
It follows from (\ref{eq:asimptotic})  that $\mathcal{H}(s)=\rho^sO(1/s) $ as $s\to\infty$.
Hence, according to the Jordan lemma \cite[chapter~V, \S2, p.439]{Lavrentev}
$$
\int\limits_{\L_{-\infty}}\mathcal{H}(s)x^{-s}ds=\int\limits_{\L_{ic}}\mathcal{H}(s) x^{-s}ds.
$$

The function $s\to\rho^sP_m(s)x^{-s}=P_m(s)e^{s\ln(\rho/x)}$ has no poles and
$$
\lim\limits_{u\to\infty}\int\limits_{-u+iu\tan\varepsilon}^{-u-iu\tan\varepsilon}\rho^s P_m(s) x^{-s}ds=0
~~\text{so that}~~\int\limits_{\L_{-\infty}}\rho^sP_m(s) x^{-s}ds=0
$$
and the definition of the $H$ function yields
\begin{equation}\label{eq:newdefinition}
H_{q,p}^{p,0}\left(x\left|\begin{array}{l}(\B,\b)\\(\A,\a)\end{array}\right.\right)
=\frac{1}{2\pi{i}}\int\limits_{\L_{ic}}\mathcal{H}(s) x^{-s}ds.
\end{equation}
According to (\ref{eq:asimptotic})
\begin{equation}\label{eq:assimHs}
\mathcal{H}(s)=\frac{\nu l_{m+1}\rho^s}{s}+\rho^s g(s),
\end{equation}
where $g(s)=O(1/s^2)$ as $s\to\infty$ in the sector $|\arg(s)|<\pi-\delta$ so that $t\to{g(c+it)}$ is absolutely
integrable on $\R$ for $c>\gamma$. Then we get the Mellin  pair
$$
v(x)=\frac{1}{2\pi{i}}\int_{c-i\infty}^{c+i\infty}(x/\rho)^{-s}g(s)ds~~\text{and}~~\int_{0}^{\infty}x^{s-1}v(x)dx=\rho^sg(s).
$$
It follows from \cite[\S12 (6)]{Boghner} that
$$
h(x)=\frac{{\nu}l_{m+1}}{2\pi{i}}\int\limits_{c-i\infty}^{c+i\infty}\frac{\rho^{s}}{s}x^{-s}ds
=\frac{{\nu}l_{m+1}}{2\pi}\int\limits_{-\infty}^{+\infty}\frac{e^{(c+it)\log(\rho/x)}}{c+it}dt=
\left\{\!\!\!\begin{array}{l} \nu l_{m+1},~~0<x<\rho,\\[10pt]
0,~~x>\rho.\end{array}\right.
$$
On the other hand, (\ref{eq:newdefinition}) implies that
$$
H_{q,p}^{p,0}\left(x\left|\begin{array}{l}(\B,\b)\\(\A,\a)\end{array}\right.\right)=h(x)+v(x)
$$
for all  $x>0$, $x\neq\rho$. Direct computation reveals that
$$
\int\limits_{0}^{\infty}x^{s-1}H_{q,p}^{p,0}\left(x\left|\begin{array}{l}
(\B,\b)\\(\A,\a)\end{array}\right.\right)dx
=\int\limits_{0}^{\rho}\nu l_{m+1}x^{s-1}dx+\rho^sg(s)=\frac{\nu l_{m+1}\rho^s}{s}+\rho^sg(s)
$$
which is equivalent to
\begin{equation}\label{eq:mellinpart}
\int_{0}^{\rho} x^{s-1}H_{q,p}^{p,0}\left(x\left|\begin{array}{l}(\B,\b)\\(\A,\a)\end{array}\right.\right)dx=\mathcal{H}(s)
~\text{for}~s>\max(0,\gamma).
\end{equation}
To remove the restriction $s>0$ note that the asymptotics of the $H$ function  given, for instance, in \cite[Theorem~1.11, p.19]{KilSaig}
implies that the left hand side of (\ref{eq:mellinpart}) is holomorphic in the half-plane $\Re{s}>\gamma$ and hence
coincides there with $\mathcal{H}(s)$.~$\hfill\square$

In our recent paper \cite{KPCMFT} we found conditions for the gamma ratio
$$
W(x)=\frac{\prod\nolimits_{i=1}^{p}\Gamma(A_ix+a_i)}{\prod\nolimits_{j=1}^{q}\Gamma(B_jx+b_j)}
$$
to be completely monotonic. According to Bernstein's theorem $W(x)$ can then be represented by the Laplace
transform of a nonnegative measure.  We showed that this measure for $\mu>0$ is given by
(\ref{eq:measureH}) and conjectured that it is given by (\ref{eq:limitmeasure}) for $\mu=0$, see
\cite[Conjecture~2]{KPCMFT}.  Immediate corollary of Theorem~~\ref{th:Fox} confirms this conjecture.
\begin{corollary}\label{cr:represent}
Suppose $\mu=0$,  $\rho\le1$ and $W(x)$ is completely monotonic. Then
$$
W(x)=\int_{[-\log(\rho),\infty)}e^{-xt}dv(t),
$$
where
$$
dv(t)={\nu}\delta_{-\log(\rho)}(t)+H_{q,p}^{p,0}\left(e^{-t}\left|\begin{array}{l}(\B,\b)\\(\A,\a)\end{array}\right.\right)\!dt.
$$
Here $\nu$  is defined in \emph{(\ref{eq:nu})} and  $\delta_{w}$ denotes the unit mass concentrated at the point $w$.
\end{corollary}

\textbf{Proof.} Application of Theorem~\ref{th:Fox} for $m=0$ yields
$$
\frac{\prod\nolimits_{k=1}^{p}\Gamma(A_ks+a_k)}{\prod\nolimits_{j=1}^{q}\Gamma(B_js+b_j)}=
\int_{0}^{\rho}H_{q,p}^{p,0}\left(x\left|\begin{array}{l}(\B,\b)\\(\A,\a)\end{array}\right.\right)x^{s-1}dx
+\nu\rho^s.
$$
It remains to apply the substitution $x=e^{-t}$ and notice that
$$
\nu\rho^s=\nu\int_{0}^{\infty}e^{-ts}d\delta_{\log(1/\rho)}(t).~~\square
$$

Next, suppose we have a family $dv_{\mathbf{c}}(x)$ of signed Borel measures supported on $[a,b]$
indexed by a real vector $\mathbf{c}$ belonging to some set $S\subset\R^n$. Assume further that $\mathbf{c}_0$
is a limit point of $S$. We recall the definition of the weak convergence \cite[Definition 8.1.1]{Bogachev}:
the signed measures $dv_{\mathbf{c}}(x)$ supported on  $[a,b]$ converge weakly to the signed measure $dv_0(x)$
as $\mathbf{c}\to\mathbf{c_0}$, $\mathbf{c}\in{S}$, if for any continuous function  $f(x)$ defined on $[a,b]$
$$
\lim\limits_{S\ni\mathbf{c}\to\mathbf{c}_0}\int_a^bf(x)dv_{\mathbf{c}}(x)=\int_a^bf(x)dv_0(x),
$$
and the limit is independent of the path.  It is relevant to remark here that the weak convergence of $dv_{\mathbf{c}}$
to the limit measure $dv_0$ on the interval $[a,b]$ does not generally imply the weak convergence of its restriction to a subinterval $[d,b]$, $a<d<b$,
to the restriction of $dv_0$ to $[d,b]$.  This is seen from the following example. The family
$$
dv_{\a}(x)=dx+\delta_{(0.5-|\a-\a_0|)}
$$
converges weakly on $[0,1]$ to the measure $dx+\delta_{0.5}$ as $\a\to\a_0$ but its restriction to $[0.5,1]$
converges weakly to the measure $dx$ on $[0.5,1]$.

In the next theorem the notation
$$
(\a^*,\b^*)=(a_1^*,...,a_p^*,b_1^*,...,b_q^*),~~~\mu^*=\sum\limits_{j=1}^{q}b_j^*-\sum\limits_{k=1}^{p}a_k^*+\frac{p-q}{2}
$$
is used.
\begin{theorem}\label{th:measureconverging}
Set $S=\{(\a,\b)\in\R^{p+q}:\a\ge0,~\mu>0\}$ and suppose $\Delta=\mu^*=0$.
Then the signed measures
$$
dv_{(\a,\b)}(x)=H_{q,p}^{p,0}\left(x\left|\begin{array}{l}(\B,\b)\\(\A,\a)\end{array}\right.\right)\!dx
$$
converge weakly on any subinterval $[\eta,\rho]$, $0\le\eta<\rho$ to the signed measure
$$
dv_{*}(x)={\nu^*\rho}\delta_\rho(x)+H_{q,p}^{p,0}\left(x\left|\begin{array}{l}(\B,\b^*)\\(\A,\a^*)\end{array}\right.\right)\!dx
$$
as $S\ni(\a,\b)\to(\a^*,\b^*)$.
Here $\rho$ is defined in \emph{(\ref{eq:rho})} and
$$
\nu^*=(2\pi)^{(p-q)/2}\prod\nolimits_{k=1}^{p}A_k^{a_k^*-1/2}\prod\nolimits_{j=1}^{q}B_j^{1/2-b_j^*}.
$$
\end{theorem}

\textbf{Proof.} Fix $p$, $q$, $\mathbf{A}$, $\mathbf{B}$. Denote
$$
H(x,\a,\b)=H_{q,p}^{p,0}\left(x\left|\begin{array}{l}(\B,\b)\\(\A,\a)\end{array}\right.\right).
$$
Write $D((\a^*,\b^*),r)$ for the closed disk of radius $r$ centered at the point $(\a^*,\b^*)$ in the space $\R^{p+q}$.
We divide the proof in three steps.

\textbf{Step~1}. We will show that $\int_0^\rho |dv_{(\a,\b)}|$ is uniformly bounded in a neighborhood of $(\a^*,\b^*)$.
The first few terms of the Stirling series \cite[p.29, 31]{ParKam} read
\begin{equation}\label{eq:paris}
\log\Gamma(z)=(z-1/2)\log z-z+\log\sqrt{2\pi}+\frac{\mathcal{B}_2}{2z}+H_2(z),
\end{equation}
where $|\arg{z}|<\pi$ and
$$
H_2(z)=\int_0^{\infty}\frac{\mathcal{B}_4-\mathcal{B}_4(t-[t])}{4(t+z)^4}dt,~~~~~\mathcal{B}_r=\mathcal{B}_r(0)~\text{are Bernoulli numbers}.
$$
Choose $\delta\in(0,1)$. Write $s=c+it$, $0<c<1$, $t\in\R$, and suppose that
$$
|s|>R=\frac{1}{\delta}\max\limits_{i,j}\left(\frac{a_i}{A_i},\frac{b_j}{B_j}\right),~~~\pi/4<|\arg{s}|<\pi/2.
$$
Standard series manipulations yield
\begin{multline}\label{eq:loggamma}
\log\Gamma(A_is+a_i)=A_is\log{s}+s\log{A_i^{A_i}}-A_is+(a_i-1/2)\log{s}+(a_i-1/2)\log{A_i}+\log\sqrt{2\pi}
\\
+\frac{1}{s}\left(\frac{-a_i^2}{2A_i}+\frac{{\cal B}_2}{2A_i}+\frac{(a_i-1/2)a_i}{A_i}\right)+
\\
\frac{1}{s^2}
\left[\sum\limits_{n=3}^\infty\frac{(-1)^{n+1}a_i^n}{n A_i^{n-1}s^{n-3}}
+\sum\limits_{n=2}^\infty\frac{(-1)^{n+1}(a_i-1/2)a_i^n}{n A_i^{n}s^{n-2}}
+\sum\limits_{n=1}^\infty\frac{(-1)^{n}{\cal B}_2a_i^n}{2A_i^{n+1} s^{n-1}}\right]+H_2(A_is+a_i).
\end{multline}
The function $\mathcal{B}_4-\mathcal{B}_4(t-[t])$ is clearly bounded on $[0,\infty)$ so that
$$
|H_2(A_is+a_i)|\leq M_0 \int\limits_0^\infty
\frac{dt}{|t+A_is+a_i|^4}=M_0\int\limits_0^\infty
\frac{dt}{|s|^4|t/|s|+A_is/|s|+a_i/|s||^4}.
$$
for some constant $M_0$ and $|s|>R$.  Change of variable $u=t/|s|+a_i/|s|$, $s/|s|=e^{i\varphi}$ gives
\begin{multline}\label{eq:rem1}
|H_2(A_is+a_i)|\leq M_0\!\!\!\int\limits_{a_i/|s|}^\infty
\frac{du}{|s|^3|u+e^{i\varphi}A_i|^4}\le
M_0\!\!\!\int\limits_{-\infty}^\infty\frac{du}{|s|^3|u+e^{i\varphi}A_i|^4}
\\[5pt]
\le\frac{M_0}{|s|^3}\biggl\{\int\nolimits_{0}^\infty
\frac{du}{\left(u^2+A_i^2\right)^2}+\int\nolimits_{-\infty}^0
\frac{du}{\left((u+ A_i/\sqrt{2})^2+A_i^2/2\right)^2}\biggr\}=\frac{M_1}{|s|^3},
\end{multline}
where we used $\pi/4<\varphi<\pi/2$. Similar relations clearly hold for $\log\Gamma(B_j s+b_j)$.  These relations
imply
\begin{equation}\label{eq:rep}
x^{-s}\frac{\prod\nolimits_{k=1}^{p}\Gamma(A_ks+a_k)}
{\prod\nolimits_{j=1}^{q}\Gamma(B_j s+b_j)}
=e^{s\log(\rho/x)}{\nu}s^{-\mu}\left(1+\frac{r(\a,\b)}{s}+\frac{R(s,\a,\b)}{s^2}\right),
\end{equation}
where  $r(\a,\b)$ is a continuous function of the vector $(\a,\b)$ and $R(s,\a,\b)$ is bounded
for $|s|\ge\delta_1$ and $(\a,\b)\in{D((\a^*,\b^*),\varepsilon)}$ with some positive $\delta_1,\varepsilon$.
On the other hand the function
$$
R(s,\a, \b)=\frac{s^{\mu+2}}{\nu \rho^s}\frac{\prod\nolimits_{k=1}^{p}\Gamma(A_ks+a_k)}
{\prod\nolimits_{j=1}^{q}\Gamma(B_j s+b_j)}-s^2- sr(\a,\b)
$$
is continuous and hence bounded on the set $\{(\a,\b)\in{D((\a^*,\b^*),\varepsilon)}, s=c+it, |s|\leq\delta_1\}$ which implies that (\ref{eq:rep}) 
holds for all $s=c+it$, $t\in\R$.  Applying (\ref{eq:rep}) to the definition (\ref{eq:Fox}) of Fox's $H$ function we get
\begin{multline}
 |H(x,\a,\b)| =
\frac{1}{2\pi}\biggl|\,\,\int\limits_{\L_{ic}}e^{s\log(\rho/x)}{\nu}
s^{-\mu}\left(1+\frac{r(\a,\b)}{s}+\frac{R(s,{\a,\b})}{s^2}\right)ds\biggr|\le
\\[5pt]
\frac{1}{2\pi}\biggl|\int\nolimits_{-\infty}^{\infty} \frac{\nu e^{(c+it)\log (\rho/x)} }{(c+it)^\mu}dt
+r({\a,\b})\int\nolimits_{-\infty}^{\infty} \frac{\nu  e^{(c+it)\log(\rho/x)} }{(c+it)^{\mu+1}}dt\biggr|
+M_2\!\int\nolimits_{-\infty}^{\infty}\frac{e^{c\log(\rho/x)}dt}{(c^2+t^2)^{1+\mu/2}}
\end{multline}
in the punctured $\varepsilon$-neighborhood of $(\a^*,\b^*)$ and for some $M_2>0$.
The next identity can be found in \cite[\S12, formula (6)]{Boghner}
$$
\dfrac{1}{2\pi}\int\limits_{-\infty}^{+\infty}\frac{e^{(c+it)\log(\rho/x)}}{(c+it)^\mu}dt=
\dfrac{1}{\Gamma(\mu)}\left(\log\frac{\rho}{x}\right)^{\mu-1},
$$
where $0<x<\rho$ and $\mu>0$.  Consequently, there exist certain constants $M_3$, $M_4$, $M_5$ such that
\begin{equation}\label{eq:estimateH}
|H(x, \a, \b)|\leq M_3\frac{\log (\rho/x))^{\mu-1}}{\Gamma(\mu)}
 + M_4\frac{\log(\rho/x))^{\mu}}{\Gamma(\mu+1)}+M_5\left(\frac{x}{\rho}\right)^{-c}
\end{equation}
in the punctured $\varepsilon$-neighborhood of the point $(\a^*,\b^*)$.
In view of $\int_0^\rho(\log(\rho/x))^{\mu-1}dx=\Gamma(\mu)$ we obtain
\begin{equation}\label{eq:estimateH1}
\int_0^\rho |dv_{(\a,\b)}(x)|\leq\int_0^\rho|H(x,\a,\b)|dx\leq M_3+M_4+M_5\rho/(1-c)=C.
\end{equation}
Using (\ref{eq:newdefinition}) for $m=0$ we get
$$
|H(x,\a^*,\b^*)|=\frac{1}{2\pi}\biggl|\,\int\limits_{\L_{ic}}e^{s\log(\rho/x)}\nu
\left(\frac{r(\a^*,\b^*)}{s}+\frac{R(s,{\a^*,\b^*})}{s^2}\right)ds\biggr|
$$
so that inequalities (\ref{eq:estimateH}) and  (\ref{eq:estimateH1}) also hold at the point  $(\a^*,\b^*)$.

\textbf{Step~2}. We now prove the theorem for $\eta=0$.  Consider arbitrary continuous function $f(x)$ defined on $[0,\rho]$.
Fix $\varepsilon_1>0$.  Weierstrass approximation theorem guarantees existence of the polynomial $P_m(x)$ such
that
$$
|P_m(x)-f(x)|<\frac{\varepsilon_1}{3C}~\text{for}~x\in[0,\rho],
$$
where $C$ is the bound from (\ref{eq:estimateH1}).  According to \cite[Theorem~7]{Fox} and Theorem~\ref{th:Fox} for $m=0$
we have
$$
\int_0^\rho x^n dv_{(\a,\b)}(x)=\frac{\prod\nolimits_{k=1}^{p}\Gamma(A_k(n+1)+a_k)}{\prod\nolimits_{j=1}^{q}\Gamma(B_j(n+1)+b_j)},
~~~\int_0^\rho x^n dv_{*}(x)=\frac{\prod\nolimits_{k=1}^{p}\Gamma(A_k(n+1)+a_k^*)}{\prod\nolimits_{j=1}^{q}\Gamma(B_j(n+1)+b_j^*)}
$$
for any $n\in\N_0$. In view of these formulas and continuity of the gamma function we have:
\begin{equation}\label{eq:lim1}
\int_0^{\rho}P_m(x)dv_{(\a,\b)}(x)=\int_0^{\rho}P_m(x)dv_{*}(x)+o(1),~~~(\a,\b)\to(\a^*,\b^*),~~~(\a,\b)\in{S}.
\end{equation}
Hence, in some punctured neighborhood of the point $(\a^*,\b^*)$ of radius $\lambda$ we have
$$
\left|\int_0^\rho  P_m(x)dv_{(\a,\b)}(x)-\int_0^\rho P_m(x)
dv_{*}(x)\right|<\frac{\varepsilon_1}{3}.
$$
Then for any $(\a,\b)\in{S}$ in the punctured neighborhood of $(\a^*,\b^*)$ of radius $\min(\lambda,\varepsilon)$  we get
\begin{multline*}
\left|\int_0^\rho f(x)dv_{(\a,\b)}(x)-\int_0^\rho
f(x)dv_{*}(x)\right|\leq \int_0^\rho \left|f(x)-P_m(x)\right||dv_{(\a,\b)}(x)|
\\
+\left|\int_0^\rho P_m(x) dv_{(\a,\b)}(x)-\int_0^{\rho}P_m(x)dv_{*}(x)\right|
\\
+\int_0^\rho\left|f(x)-P_m(x)\right||dv_{*}(x)|\leq
C\frac{\varepsilon_1}{3C}+\frac{\varepsilon_1}{3}+C\frac{\varepsilon_1}{3C}=\varepsilon_1,
\end{multline*}
which proves the weak convergence of $dv_{(\a,\b)}(x)$ on $[0,\rho]$.

\textbf{Step~3}. Suppose now that $0<\eta<\rho$. For any continuous function $g(x)$ defined on $[\eta, \rho]$
we can build the function
$$
g_\sigma(x)=\left\{\begin{array}{l} 0,\ 0\leq x\leq \eta-\sigma,
\\
g(\eta)((x-\eta)/\sigma+1),\ \eta-\sigma< x<\eta,
\\
g(x),\ \eta\leq x\leq \rho,
\end{array}\right.
$$
where $\sigma>0$ is sufficiently small. It follows from (\ref{eq:estimateH}) that the function $x\to{H(x,\a,\b)}$ is bounded
on $[\eta-\sigma,\eta]$ for $(\a,\b)$ in some neighborhood of  $(\a^*,\b^*)$. Let us call this bound $M$.  Then
\begin{multline*}
\left|\int_\eta^\rho g(x)dv_{(\a,\b)}(x)-\int_\eta^\rho
g(x)dv_{*}(x)\right|\leq \left|\int_0^\rho
g_\sigma(x)dv_{(\a,\b)}(x)-\int_0^\rho g_\sigma(x)dv_{*}(x)\right|
\\
+\int_{\eta-\sigma}^\eta |g_\sigma(x)||H(x,\a,\b)|dx+\int_{\eta-\sigma}^\eta|g_\sigma(x)||H(x,\a^*,\b^*)|dx
\\
\leq\left|\int_0^\rho g_\sigma(x)dv_{(\a,\b)}(x)-\int_0^\rho
g_\sigma(x) dv_{*}(x)\right|+ 4M|g(\eta)|\sigma.
\end{multline*}
In view of the weak convergence of $dv_{(\a,\b)}$ on $[0,\rho]$ we obtain
$$
\limsup\limits_{S\ni(\a,\b)\to (\a^*,\b^*)}\left|\int_\eta^{\rho}g(x)dv_{(\a,\b)}(x)
-\int_\eta^{\rho}g(x)dv_{*}(x)\right|\leq4M|g(\eta)|\sigma.
$$
As $\sigma$ is arbitrary,
$$
\lim\limits_{S\ni(\a,\b)\to(\a^*,\b^*)}\left|\int_\eta^\rho g(x)dv_{(\a,\b)}(x)
-\int_\eta^{\rho}g(x)dv_{*}(x)\right|=0.~~~~~\square
$$

The above theorem leads to an extension of the integral equation for the delta neutral $H$ function obtained in
\cite[Theorem~8]{KPCMFT} to the case $\mu=0$.

\begin{corollary}\label{cr:intequation}
Let $\A,\B>0$, $\a,\b\ge0$, $\mu=0$, $\rho\le1$ and
$\sum_{i=1}^{p}A_i=\sum_{j=1}^{q}B_j$. Then the identity
\begin{equation}\label{eq:Hid}
H_{q,p}^{p,0}\left(x\left|\begin{array}{l}(\B,\b)\\(\A,\a)\end{array}\right.\right)=\frac{1}{\log(\rho/x)}\left(\nu
Q(x/\rho)+\int_{x}^{\rho}
H_{q,p}^{p,0}\left(u\left|\begin{array}{l}(\B,\b)\\(\A,\a)\end{array}\right.\right)\frac{Q(x/u)}{u}du\right),
\end{equation}
holds for all $x\in(0,\rho)$ with $Q$ defined in \emph{(\ref{eq:Qu-def})}.
\end{corollary}

\textbf{Proof.}  Let $\varepsilon>0$, $\b_\varepsilon=(b_1+\varepsilon, b_2,....,b_q)$.
Then $\mu_\varepsilon=\mu+\varepsilon>0$ so that the $H$ function with parameters $\A,\B$, $\a,\b_\varepsilon$  satisfies
the integral equation \cite[(25)]{KPCMFT}
\begin{equation}\label{H-ident}
\log(\rho/x)H_{q,p}^{p,0}\left(x\left|\begin{array}{l}(\B,\b_\varepsilon)\\(\A,\a)\end{array}\right.\right)=\int_{x}^{\rho}
H_{q,p}^{p,0}\left(u\left|\begin{array}{l}(\B,\b_\varepsilon)\\(\A,\a)\end{array}\right.\right)\frac{Q_\varepsilon(x/u)}{u}du,
\end{equation}
where
\begin{equation}\label{eq:Q2}
Q_\varepsilon(t)=\sum\limits_{i=1}^{p}\frac{t^{a_i/A_i}}{1-t^{1/A_i}}-\frac{t^{(b_1+\varepsilon)/B_1}}{1-t^{1/B_1}}-
\sum\limits_{j=2}^{q}\frac{t^{b_j/B_j}}{1-t^{1/B_j}}.
\end{equation}
It is tempting to apply Theorem~\ref{th:measureconverging} to the integral on the right of (\ref{H-ident}) to arrive at
the required conclusion immediately. However, the dependence of $Q_\varepsilon(x/u)$  on $\varepsilon$ is not permitted by
conditions of this theorem, so that one more step is needed.  Denote
$$
H_\varepsilon(x)=H_{q,p}^{p,0}\left(x\left|\begin{array}{l}(\B,\b_\varepsilon)\\(\A,\a)\end{array}\right.\right),\
H(x)=H_{q,p}^{p,0}\left(x\left|\begin{array}{l}(\B,\b)\\(\A,\a)\end{array}\right.\right).
$$
The triangle inequality yields
\begin{multline}\label{eq:treyg1}
\left|\int_{x}^{\rho} H_\varepsilon(u)\frac{Q_\varepsilon
(x/u)}{u}du-\int_{x}^{\rho} H(u)\frac{Q(x/u)}{u}du-\nu
Q(x/\rho)\right|\leq
\\
\int_{x}^{\rho}
\left|\frac{H_\varepsilon(u)}{u}\right||Q_\varepsilon(x/u)-Q(x/u)|du+
\left|\int_{x}^{\rho}(H_\varepsilon(u)-H(u))\frac{Q(x/u)}{u}du-\nu Q(x/\rho)\right|.
\end{multline}
We have
$$
Q_\varepsilon(t)-Q(t)=\frac{t^{b_1/B_1}(t^{\varepsilon/B_1}-1)}{1-t^{1/B_1}}
$$
so that
$$
\sup\limits_{t\in[x/\rho,1]}|Q_\varepsilon(t)-Q(t)|\leq
\sup\limits_{t\in[x/\rho,1]}
\frac{1-t^{\varepsilon/B_1}}{1-t^{1/B_1}}=\frac{1-(x/\rho)^{\varepsilon/B_1}}{1-(x/\rho)^{1/B_1}}
$$
and hence
$$
\int_{x}^{\rho}
\left|\frac{H_\varepsilon(u)}{u}\right||Q_\varepsilon(x/u)-Q(x/u)|du\leq
\frac{1-(x/\rho)^{\varepsilon/B_1}}{1-(x/\rho)^{1/B_1}}\int_{x}^{\rho}
\frac{|H_\varepsilon(u)|}{x} du.
$$
The last estimate combined with (\ref{eq:estimateH1}) leads to the conclusion that the first term on the right hand side of (\ref{eq:treyg1})
tends to zero as $\varepsilon\to0$. The second term also goes to zero as seen by integrating the difference $H_\varepsilon(u)-H(u)$
term by term and applying Theorem~\ref{th:measureconverging} to the first term.
The claim now follows on limit transition $\varepsilon\to0$ in (\ref{H-ident}). $\hfill\square$

\paragraph{3. New identities for Meijer's $G$ function.}
Meijer's $G$ function is a particular case of Fox's $H$ function when $\A=\B=(1,\ldots,1)$.  The condition $\Delta=0$
then forces the assumption $p=q$ for $G$ to be delta neutral.
It also  follows from (\ref{eq:rho}) that $\rho=1$.
Under the conditions $\a,\b\ge0$, $0<x<1$, equation (\ref{eq:ident1}) takes the following form
\begin{equation}\label{eq:G-represent1}
\log(1/x)G^{p,0}_{p,p}\left(x\,\,\vline\begin{array}{c}\!\b\\\!\a\end{array}\!\!\right)
=\int_x^{1}G^{p,0}_{p,p}\left(t\,\,\vline\begin{array}{c}\!\b\\\!\a\end{array}\!\!\right)
\sum_{k=1}^{p}\left(\frac{x^{a_k}}{t^{a_k}}-\frac{x^{b_k}}{t^{b_k}}\right)\frac{dt}{t-x}
\end{equation}
where $\mu>0$, while equation (\ref{eq:Hid}) reduces to
\begin{equation}\label{eq:G-represent2}
\log(1/x)G^{p,0}_{p,p}\left(x\,\,\vline\begin{array}{c}\!\b\\\!\a\end{array}\!\!\right)
=\sum_{k=1}^{p}\frac{\left(x^{a_k}-x^{b_k}\right)}{1-x}+\int_x^{1}
G^{p,0}_{p,p}\left(t\,\,\vline\begin{array}{c}\!\b\\\!\a\end{array}\!\!\right)
\sum_{k=1}^{p}\left(\frac{x^{a_k}}{t^{a_k}}-\frac{x^{b_k}}{t^{b_k}}\right)\frac{dt}{t-x}
\end{equation}
for $\mu=0$. Formula (\ref{eq:G-represent1}) is
nontrivial even for  $p=1$ and $p=2$ (as well as (\ref{eq:G-represent2}) for $p=2$).
Using the representations \cite[8.4.2.3, 8.4.49.22]{PBM3}
\begin{equation}\label{eq:G1011}
G^{1,0}_{1,1}\!\left(\!z~\vline\begin{array}{l}b\\a\end{array}\!\!\right)=
\frac{z^a(1-z)^{b-a-1}}{\Gamma(b-a)}
\end{equation}
and
\begin{equation}\label{eq:G2022}
G^{2,0}_{2,2}\!\left(\!z~\vline\begin{array}{l}b_1,b_2
\\a_1,a_2\end{array}\!\!\right)=\frac{z^{a_2}(1-z)^{b_1+b_2-a_1-a_2-1}}{\Gamma(b_1+b_2-a_1-a_2)}
{_2F_{1}}\!\left(\begin{array}{l}b_1-a_1,b_2-a_1
\\b_1+b_2-a_1-a_2\end{array}\!\!;1-z\right)
\end{equation}
valid for $|z|<1$, we get the following curious integral evaluations:
$$
\int_{x}^{1}t^a(1-t)^{b-a-1}\left[\frac{x^a}{t^a}-\frac{x^b}{t^b}\right]\frac{dt}{t-x}=x^{a}(1-x)^{b-a-1}\log(1/x),~~b>a,
$$
and
\begin{multline*}
\int_{x}^{1}t^{a_1}(1-t)^{\psi-1}
{_2F_{1}}\!\left(\left.\begin{array}{l}b_1-a_2,b_2-a_2\\\psi\end{array}\!\!\right|1-t\right)
\left[\frac{x^{a_1}}{t^{a_1}}+\frac{x^{a_2}}{t^{a_2}}-\frac{x^{b_1}}{t^{b_1}}-\frac{x^{b_2}}{t^{b_2}}\right]\frac{dt}{t-x}
\\
=x^{a_1}(1-x)^{\psi-1}\log(1/x){_2F_{1}}\!\left(\left.\begin{array}{l}b_1-a_2,b_2-a_2\\\psi\end{array}\!\!\right|1-x\right),~~\mu=b_1+b_2-a_1-a_2>0.
\end{multline*}
The order of summation and integration in (\ref{eq:G-represent1}) and (\ref{eq:G-represent2}) cannot be interchanged
since the resulting integrals diverge. Nevertheless, these divergent integrals can be regularized which will be done in
Theorem~\ref{th:G-integral} below.
Denote
$$
\widehat{G}^{p+1,0}_{p+1,p+1}\left(x\,\,\vline\begin{array}{c}\!b_0,\b\\\!u,\a\end{array}\!\!\right)
=\frac{d}{du}G^{p+1,0}_{p+1,p+1}\left(x\,\,\vline\begin{array}{c}\!b_0,\b\\\!u,\a\end{array}\!\!\right).
$$
Differentiating the definition of the $G$ function we obtain:
\begin{equation}\label{eq:D1G-integral}
\widehat{G}^{p+1,0}_{p+1,p+1}\left(x\,\,\vline\begin{array}{c}\!b_0,\b
\\\!u,\a\end{array}\!\!\right)=\frac{1}{2\pi{i}}
\int\limits_{\L_{-\infty}}\!\!\frac{\Gamma'(u\!+\!s)\Gamma(a_1\!+\!s)\dots\Gamma(a_p\!+\!s)}
{\Gamma(b_{1}\!+\!s)\dots\Gamma(b_{p+1}\!+\!s)}x^{-s}ds.
\end{equation}

\begin{theorem}\label{th:G-integral}
Suppose $\a,\b\ge0$ and $0<x<1$. Then the $G$ function
satisfies the following equations\emph{:}
\begin{gather}
\log(1/x)G^{p,0}_{p,p}\left(x\,\,\vline\begin{array}{c}\!\b\\\!\a\end{array}\!\!\right)
=\sum_{k=1}^{p}\left\{\widehat{G}^{p+1,0}_{p+1,p+1}\left(x\,\,\vline\begin{array}{c}\!b_k,\b\\\!b_k,\a\end{array}\!\!\right)
-\widehat{G}^{p+1,0}_{p+1,p+1}\left(x\,\,\vline\begin{array}{c}\!a_k,\b\\\!a_k,\a\end{array}\!\!\right)\right\}\label{eq:G-represent}
\\[8pt]
=\frac{1}{2\pi{i}}\int\limits_{\L_{-\infty}}\!\!\frac{\Gamma(a_1\!+\!s)\dots\Gamma(a_p\!+\!s)}
{\Gamma(b_{1}\!+\!s)\dots\Gamma(b_{p}\!+\!s)}\left(\sum_{k=1}^{p}(\psi(b_k+s)-\psi(a_k+s)\right)x^{-s}ds\label{eq:psi-integral}
\end{gather}
if $\mu>0$ and
\begin{multline}\label{eq:Gzero}
\log (1/x)
G^{p,0}_{p,p}\left(x\,\,\vline\begin{array}{c}\!\b\\\!\a\end{array}\!\!\right)
\!=\!\sum_{k=1}^{p}\frac{\left(x^{a_k}-x^{b_k}\right)}{1-x}+
\sum_{k=1}^{p}\left\{\widehat{G}^{p+1,0}_{p+1,p+1}\left(x\,\,\vline\begin{array}{c}\!b_k,\b\\\!b_k,\a\end{array}\!\!\right)
-\widehat{G}^{p+1,0}_{p+1,p+1}\left(x\,\,\vline\begin{array}{c}\!a_k,\b\\\!a_k,\a\end{array}\!\!\right)\right\}
\\[8pt]
=\sum_{k=1}^{p}\frac{\left(x^{a_k}-x^{b_k}\right)}{1-x}+\frac{1}{2\pi{i}}\int\limits_{\L_{-\infty}}\!\!\frac{\Gamma(a_1\!+\!s)\dots\Gamma(a_p\!+\!s)}
{\Gamma(b_{1}\!+\!s)\dots\Gamma(b_{p}\!+\!s)}\left(\sum_{k=1}^{p}(\psi(b_k+s)-\psi(a_k+s)\right)x^{-s}ds
\end{multline}
if $\mu=0$.  Here $\psi(z)=\Gamma'(z)/\Gamma(z)$ denotes the digamma function.
\end{theorem}

\textbf{Proof.} Write the right hand side of (\ref{eq:G-represent1}) as
\begin{equation}\label{eq:limitG}
\lim\limits_{h\to 0}
\sum_{k=1}^{p}\left\{x^{a_k}\int_x^{1}G^{p,0}_{p,p}\left(t\,\,\vline\begin{array}{c}\!\b\\\!\a\end{array}\!\!\right)(t-x)^{h-1}t^{-a_k}
dt-x^{b_k}\int_x^{1}G^{p,0}_{p,p}\left(t\,\,\vline\begin{array}{c}\!\b\\\!\a\end{array}\!\!\right)(t-x)^{h-1}t^{-b_k}
dt\right\}.
\end{equation}
According to \cite[Formula~2.24.3, p.~293]{PBM3} we get
$$
x^{a_k}\int_x^{1}G^{p,0}_{p,p}\left(t\,\,\vline\begin{array}{c}\!\b\\\!\a\end{array}\!\!\right)(t-x)^{h-1}t^{(1-a_k)-1}
dt=x^h \Gamma(h)
G^{p+1,0}_{p+1,p+1}\left(x\,\,\vline\begin{array}{c}\!a_k,\b
\\\!a_k-h,\a\end{array}\!\!\right).
$$
and
$$
x^{b_k}\int_x^{1}G^{p,0}_{p,p}\left(t\,\,\vline\begin{array}{c}\!\b\\\!\a\end{array}\!\!\right)(t-x)^{h-1}t^{(1-b_k)-1}
dt=x^h \Gamma(h)
G^{p+1,0}_{p+1,p+1}\left(x\,\,\vline\begin{array}{c}\!b_k,\b
\\\!b_k-h,\a\end{array}\!\!\right).
$$
By Taylor's theorem
$$
G^{p+1,0}_{p+1,p+1}\left(x\,\,\vline\begin{array}{c}\!a_k,\b
\\\!a_k-h,\a\end{array}\!\!\right)= G^{p,0}_{p,p}\left(x\,\,\vline\begin{array}{c}\!\b
\\\!\a\end{array}\!\!\right)-h\widehat{G}^{p+1,0}_{p+1,p+1}\left(x\,\,\vline\begin{array}{c}\!a_k,\b
\\\!a_k,\a\end{array}\!\!\right)+o(h)~\text{as}~h\to0
$$
and
$$
G^{p+1,0}_{p+1,p+1}\left(x\,\,\vline\begin{array}{c}\!b_k,\b
\\\!b_k-h,\a\end{array}\!\!\right)= G^{p,0}_{p,p}\left(x\,\,\vline\begin{array}{c}\!\b
\\\!\a\end{array}\!\!\right)-h\widehat{G}^{p+1,0}_{p+1,p+1}\left(x\,\,\vline\begin{array}{c}\!b_k,\b
\\\!b_k,\a\end{array}\!\!\right)+o(h)~\text{as}~h\to0.
$$
As $\lim_{h\to0}h\Gamma(h)=1$ we conclude that the limit (\ref{eq:limitG}) equals
$$
\sum_{k=1}^{p}\left\{\widehat{G}^{p+1,0}_{p+1,p+1}\left(x\,\,\vline\begin{array}{c}\!b_k,\b
\\\!b_k,\a\end{array}\!\!\right)-\widehat{G}^{p+1,0}_{p+1,p+1}\left(x\,\,\vline\begin{array}{c}\!a_k,\b
\\\!a_k,\a\end{array}\!\!\right)\right\}.
$$
Formula (\ref{eq:psi-integral}) now follows on substituting (\ref{eq:D1G-integral}) into (\ref{eq:G-represent}).
The proof of (\ref{eq:Gzero}) is similar. $\hfill\square$

\textbf{Remark.} N{\o}rlund's identity (\ref{eq:Norlund1}) for the $G$ function in (\ref{eq:G-represent}) takes the form
\begin{equation*}\label{eq:Gp+1}
G^{p+1,0}_{p+1,p+1}\!\left(\!x~\vline\begin{array}{l}b_k,\b\\u,\a\end{array}\!\!\right)=x^{u}(1-x)^{\mu+b_k-u-1}
\sum\limits_{n=0}^{\infty}\frac{g_n(\a;b_k,\b)}{\Gamma(\mu+b_k-u+n)}(1-x)^n.
\end{equation*}
The coefficients $g_n(\a;b_k,\b)$ are given in (\ref{eq:Norlund-explicit}). Further details regarding
these coefficients can be found in \cite{Norlund} and our recent paper \cite{KPNorlund}.
Computing the derivative in $u$ and setting $u=b_k$ we obtain
\begin{multline*}
\widehat{G}^{p+1,0}_{p+1,p+1}\!\left(\!x~\vline\begin{array}{l}b_k,\b\\b_k,\a\end{array}\!\!\right)
=(\log(x)-\log(1-x))x^{b_k}(1-x)^{\mu-1}\sum\limits_{n=0}^{\infty}\frac{g_n(\a;b_k,\b)}{\Gamma(\mu+n)}(1-x)^n
\\
+x^{b_k}(1-x)^{\mu-1}\sum\limits_{n=0}^{\infty}\frac{g_n(\a;b_k,\b)\psi(\mu+n)}{\Gamma(\mu+n)}(1-x)^n
\end{multline*}
and (setting $u=a_k$)
\begin{multline*}
\widehat{G}^{p+1,0}_{p+1,p+1}\!\left(\!x~\vline\begin{array}{l}a_k,\b\\a_k,\a\end{array}\!\!\right)
=(\log(x)-\log(1-x))x^{a_k}(1-x)^{\mu-1}\sum\limits_{n=0}^{\infty}\frac{g_n(\a;a_k,\b)}{\Gamma(\mu+n)}(1-x)^n
\\
+x^{a_k}(1-x)^{\mu-1}\sum\limits_{n=0}^{\infty}\frac{g_n(\a;a_k,\b)\psi(\mu+n)}{\Gamma(\mu+n)}(1-x)^n.
\end{multline*}
It is easily seen on setting $u=b_k$ or $u=a_k$ that
$$
x^{b_k}(1-x)^{\mu-1}\sum\limits_{n=0}^{\infty}\frac{g_n(\a;b_k,\b)}{\Gamma(\mu+n)}(1-x)^n
=x^{a_k}(1-x)^{\mu-1}\sum\limits_{n=0}^{\infty}\frac{g_n(\a;a_k,\b)}{\Gamma(\mu+n)}(1-x)^n
=G^{p,0}_{p,p}\!\left(\!x~\vline\begin{array}{l}\b\\\a\end{array}\!\!\right).
$$
Hence, (\ref{eq:G-represent}) takes the form
$$
\log(1/x)G^{p,0}_{p,p}\left(x\,\,\vline\begin{array}{c}\!\b\\\!\a\end{array}\!\!\right)
\!=\!(1-x)^{\mu-1}\sum\limits_{n=0}^{\infty}\frac{\psi(\mu+n)}{\Gamma(\mu+n)}(1-x)^n
\sum_{k=1}^{p}\{x^{b_k}g_n(\a;b_k,\b)-x^{a_k}g_n(\a;a_k,\b)\}.
$$
In particular, for $p=2$ we get using (\ref{eq:G2022}) and \cite[(2.10)]{Norlund}:
\begin{multline*}
x^{a_1}\log(1/x){_2F_{1}}\!\left(\!\!\left.\begin{array}{l}b_1-a_2,b_2-a_2\\\mu\end{array}\!\!\right|1-x\right)
\\
=\!\sum\limits_{n=0}^{\infty}\frac{\psi(\mu+n)}{n!}\left\{x^{b_1}(\mu)_n{_3F_{2}}\!\left(\!\!\left.\begin{array}{l}-n,b_2-a_1,b_2-a_2\\\mu,\mu\end{array}\!\!\right|\!1\!\right)
\!+\!x^{b_2}(\mu)_n{_3F_{2}}\!\left(\!\!\left.\begin{array}{l}-n,b_2-a_1,b_2-a_2\\\mu,\mu\end{array}\!\!\right|\!1\!\right)\right.
\\
\left.-x^{a_1}\frac{(b_1-a_2)_n(b_2-a_2)_n}{(\mu)_n}-x^{a_2}\frac{(b_1-a_1)_n(b_2-a_1)_n}{(\mu)_n}\right\}(1-x)^n.
\end{multline*}
The ultimate  two identities are believed to be new.

\paragraph{4. Acknowledgements.}  This work has been supported by the Russian Science Foundation under project 14-11-00022.

\end{document}